\documentclass[reqno,12pt]{amsart}

\usepackage{amsmath}
\usepackage{amsfonts}
\usepackage{amssymb}
\usepackage{eufrak}
\usepackage{amscd}
\usepackage{amsthm}
\usepackage{epsfig}
\usepackage{amstext}
\usepackage[all]{xy}

\theoremstyle{plain}

 \newtheorem{theorem}{Theorem}

 \newtheorem{cor}[theorem]{Corollary}
 \newtheorem{lemma} [theorem] {Lemma}
 \newtheorem{remark}[theorem]{Remark}

\begin{document}

 \title{Lifting algebraic contractions in C*-algebras}
\author{Terry Loring and Tatiana Shulman}

\address{Department of Mathematics and Statistics, University of New Mexico,
Albuquerque, NM 87131, USA.}

\address{Department of Mathematics, Siena College, 515 Loudon Road, Loudonville, NY 12211, USA.}


\subjclass[2000]{46 L05; 46L35}

\keywords {Projective and semiprojective $C^*$-algebras, stable
relation, lifting problem}

\date{\today}

\thanks{This work was partially supported by a grant from the
Simons Foundation (208723 to Loring) and by the NordForsk Research
Network "Operator Algebras and Dynamics" (grant 11580) .}

\maketitle

\begin{abstract}Let $p$ be a polynomial in one variable. It is shown that the universal $C^*$-algebra of the relation $p(x)=0$, $\|x\| \le C$ is
semiprojective, residually finite-dimensional and has trivial
extension group.
\end{abstract}

\section*{Introduction}
Notions of projectivity and semiprojectivity for $C^*$-algebras was
introduced by Effros and Kaminker \cite{EffrKam} and, in its modern
form, by Blackadar \cite{Blackadar} as noncommutative analogues of
absolute retract and absolute neighborhood retract in topology.

A $C^*$-algebra $D$ is  {\it projective} if for any $C^*$-algebra
$A$, its ideal $I$ and every $\ast$-homomorphism $\phi:D \to A/I$,
there exists  a $\ast$-homomorphism $\tilde \phi$ such that the
diagram
$$\xymatrix{ & A \ar[d] \\ D\ar[ru]^{\tilde \phi}\ar[r]^{\phi} & A/I}$$ commutes.

A $C^*$-algebra $D$ is  {\it semiprojective}  if for any
$C^*$-algebra $A$, any increasing chain of ideals $I_1\subseteq
I_2\subseteq \ldots\;\;$ in $A$ and for every $\ast$-homomorphism
$\phi:D \to A/\;\overline{\bigcup_k I_k}$,
 there exist
 $n$ and  a $\ast$-homomorphism $\tilde \phi : D\to A/I_n$ such that the diagram
$$\xymatrix{ & A/I_n \ar[d] \\ D\ar[ru]^{\tilde \phi}\ar[r]^{\phi} &
A/\;\overline{\bigcup_k I_k}}$$ commutes.

A notion of weak semiprojectiivty was introduced by Eilers and
Loring \cite{SorenTerry}.

A $C^*$-algebra $D$ is  {\it weakly semiprojective}  if for any
sequence  of $C^*$-algebras $A_i$ and any $\ast$-homomorphism
$\phi:D \to \prod A_i /\bigoplus A_i$, there exists a
$\ast$-homomorphism  $\tilde \phi: D \to \prod A_i$ such that the
diagram
$$\xymatrix{ & \prod A_i \ar[d] \\ D\ar[ru]^{\tilde \phi}\ar[r]^{\phi} &\prod A_i /\bigoplus A_i}$$ commutes.

Examples and basic properties of projective and (weakly)
semiprojective $C^*$-algebras can be found in \cite{book}.

The notion of (weak) (semi)projectivity provides an algebraic
setting of lifting and perturbation problems for relations in
$C^*$-algebras. Namely a relation is liftable (which means that in
any quotient $C^*$-algebra $A/I$ any elements satisfying the
relation have preimages in $A$ also satisfying the relation) if and
only if the universal $C^*$-algebra of the relation is projective.

Similarly a relation is stable under small perturbations if and only
if the universal $C^*$-algebra of the relation is weakly
semiprojective.

 The only problem is that
not all relations have the universal $C^*$-algebras. However if
relations are noncommutative $\ast$-polynomial equations combined
with norm restrictions on generators $$\|x_i\|\le c_i$$ then such
system of relations defines the universal $C^*$-algebra. That is why
in (weak) (semi)projectivity questions it is important to solve
lifting problems for relations combined with norm restrictions on
generators.

It is not easy for a polynomial relation to be liftable. For
example, let us consider a polynomial in one variable. Suppose in a
$C^*$-quotient we have an element $x$ satisfying $p(x)=0$. If $p$
has a non-zero root, the spectral idempotent of $x$ on this root
also belongs to the quotient. But idempotents are not liftable:
consider for example the unit of $\mathbb C = C_0(0, 1]/C_0(0, 1)$.
It has no idempotent lift because $C_0(0, 1]$ contains no non-zero
idempotents. Thus only monomials $x^n=0$ have a chance to be
liftable. And they indeed are liftable by deep result of Olsen and
Pedersen \cite{OlsenPedersen}. Immediately there arises a question
(\cite{book}) if the universal $C^*$-algebra of $x^n=0$, $\|x\| \le
C$ is projective. It was answered positively in \cite{nilpotents}
and in this paper we will give a short proof of that (Corollary
\ref{nilp}).

As to non-monomial relations $p(x)=0$, they are known to be stable
under small perturbations \cite{Don}. Moreover the proof in
\cite{Don} can be generalized to show that relations $p(x)=0$ are
liftable from quotients of the form $A/\;\overline{\bigcup I_n}$
arising in definition of semiprojectivity,  and one is led to ask if
the universal $C^*$-algebra of $p(x)=0$, $\|x\| \le C$ is
semiprojective. In \cite{multiplicity1} it was proved in the case
when all roots of the polynomial have multiplicity more than 1. In
this paper we prove it for arbitrary polynomial (Theorem
\ref{semiproj}).

We also show that these universal $C^*$-algebras are RFD (Theorem
\ref{RFD}) and that $\ast$-homomorphisms from these $C^*$-algebras
to Calkin algebra lift to $\ast$-homomorphisms to $B(H)$ (Theorem
\ref{Calkin}). The last result is generalization of Olsen's
structure theorem for polynomially compact operators
\cite{OlsenAlg}.

Our main technical tool is a generalized spectral radius formula we
introduced in \cite{OlsenQuestion} in connection with question of
Olsen about best approximation of operators by compacts. It turns
out to be useful tool also for lifting polynomial relations combined
with restrictions on norms of generators.

\section*{A generalized spectral radius formula}

For $x\in A$, we denote by $\dot x$ its image in $A/I$ and by
$\rho(x)$ its spectral radius.

The following theorem is a generalization of spectral radius formula
of Murphy and West \cite{MurphyWest}. The spectral radius formula is
particular case of the generalized spectral radius formula when
$I=A$.

  \begin{theorem}(\cite{OlsenQuestion}) Let $A$ be a $C^*$-algebra, $I$ its ideal, $x\in A$.
  Then $$\max\{\rho(x), \|\dot x\|\} = \inf \|(1+i)x(1+i)^{-1}\|$$ (here  $\inf$ is taken over
  all $i\in I$ such that $1+i$ is invertible). If $\|\dot x\| > \rho(x)$ then the infimum in the right-hand side is attained.
  \end{theorem}

\begin{lemma}\label{basic} Let $p$ be a polynomial in one variable and $t_1, \ldots, t_k$ its roots. Let $A$ be a $C^*$-algebra, I its ideal, $x\in A/I$,
$p(x)=0$ and $\|x\|  > \max\{t_i\}$. Suppose $x$ has a lift $X\in A$
such that $p(X)=0$. Then there is a lift $\tilde X\in A$ of $x$ such
that $p(\tilde X)=0$ and $\|\tilde X\| =\|x\|.$
\end{lemma}
\begin{proof} Since $p(X)=0$, $$\rho(X) = \max\{t_i\} <
\|x\|.$$ By the generalized spectral radius formula there exists
$i\in I$ such that $$\|(1+i)X(1+i)^{-1}\| = \|x\|.$$ Let $\tilde X =
(1+i)X(1+i)^{-1}$. Then $\tilde X$ is a lift of $x$, $\|\tilde X \|
= \|x\|$ and $$p(\tilde X) = p((1+i)X(1+i)^{-1}) =
(1+i)p(X)(1+i)^{-1} = 0.$$
\end{proof}

\begin{cor}\label{nilp}(\cite{nilpotents}) The universal $C^*$-algebra $$C^*\langle x\;|\; x^n = 0, \|x\| \le C\rangle$$ is projective.
\end{cor}
\begin{proof} For $C=0$ the statement obviously holds. So let $C>0$ and let $x\in A/I$, $x^n~=~0, \; \|x\| \le
C.$ We need to show that there a lift of $x$ with the same
properties. If $x=0$ then it is obvious, so let us assume $x\neq 0$.
By \cite{OlsenPedersen} there is a lift $X$ of $x$ such that
$X^n=0$. By Lemma \ref{basic} there is a lift $\tilde X$ of $x$ such
that $(\tilde X)^n=0$ and $\|\tilde X\| = \|x\| \le C.$
\end{proof}

\section*{Semiprojectivity of the universal $C^*$-algebra $C^*\langle x\;|\; p(x) = 0, \|x\| \le
C\rangle$}

 \begin{lemma}\label{1} Let $T\in B(H)$ and $(T-t_N)^{k_N}(T-t_{N-1})^{k_{N-1}}\ldots (T-t_1)^{k_1}=0$.
 Let \begin{multline*} \\ H_1 = \ker (T-t_1) \\ H_2 = \ker (T-t_1)^2 \ominus H_1 \\ \ldots \\ H_{k_1} = \ker (T-t_1)^{k_1}\ominus H_{k_1-1} \\ H_{k_1 + 1} =
\ker (T-t_2)(T-t_1)^{k_1} \ominus H_{k_1}\\ \ldots \\
H_{k_1+\ldots + k_N} = \ker (T-t_N)^{k_N
-1}(T-t_{N-1})^{k_{N-1}}\ldots (T-t_1)^{k_1} \ominus H_{k_1+\ldots +
k_N -1}.\end{multline*}

\medskip

\noindent Then with respect to the decomposition $H = H_1 \oplus
\ldots \oplus H_{k_1+\ldots + k_N}$ the operator $T$ has
uppertriangular form with $t_1 1, \ldots, t_1 1, \ldots, t_N 1,
\ldots, t_N 1$ on the diagonal, where each $t_i 1$ is repeated $k_i$
times.
\end{lemma}
\begin{proof} If $x\in H_1$, then $Tx = t_1x$. If $x\in H_2$, then $Tx = (T-t_1)x + t_1x$, where $(T-t_1)x \in H_1$.
And so on.
\end{proof}

\medskip

\begin{lemma}\label{idempotent} Let $B\subseteq B(H)$ be a $C^*$-algebra, $b\in B$ an idempotent. Then the projection onto the range of b also belongs to $B$.
\end{lemma}
\begin{proof} By Lemma \ref{1}, $b$ can be written as \begin{equation}b = \left(\begin{array}{cc} 1 & X \\ 0 & 0 \end{array}\right).\end{equation} Hence $$\left(\begin{array}{cc} 1 + XX^* & 0 \\ 0 & 0 \end{array}\right) =  \left(\begin{array}{cc} 1 & X \\ 0 & 0 \end{array}\right) \left(\begin{array}{cc} 1 & 0 \\ X^* & 0 \end{array}\right) = bb^*  \in B.$$ Let $f$ be a continuous function on $\mathbb R$ which vanishes at $0$ and equal $1$ at $[1, \infty)$. Then
$$ \left(\begin{array}{cc} 1 & 0 \\ 0 & 0 \end{array}\right) = f\left(\left(\begin{array}{cc} 1 + XX^*& 0 \\ 0 & 0 \end{array}\right)\right) \in B.$$ From (1) it is seen that  $\left(\begin{array}{cc} 1 & 0 \\ 0 & 0 \end{array}\right)$ is exactly the projection onto the range of $b$.
\end{proof}

\medskip

Let $L_1$ be direct sum of first $k_1$ summands in the decomposition
$H = H_1 \oplus \ldots \oplus H_{k_1+\ldots + k_N -1}$ (that is
$\ker (T-t_1)^{k_1}$), $L_2$ be direct sum of next ${k_2}$ summands,
and so on. For any $m$, let  $M_m = H \ominus (L_1\oplus \ldots
\oplus L_m).$

\medskip

\begin{cor}\label{belong} Let $T\in B(H)$,   $$(T-t_N)^{k_N}(T-t_{N-1})^{k_{N-1}}\ldots (T-t_1)^{k_1}=0 $$ and subspaces $L_i$ be  as above. Then projections onto $L_i$ belong to $C^*(T, 1)$.
\end{cor}
\begin{proof} By transposing factors in the product $(T-t_N)^{k_N}(T-t_{N-1})^{k_{N-1}}\ldots (T-t_1)^{k_1}$  the general case can be reduced to the case $i=1$. So let us prove that the projection onto $L_1$ belongs to $C^*(T, 1)$.  Since $t_1$ is an isolated point of $\sigma(T)$, there exists   the spectral idempotent $Q$  corresponding to $t_1$, that is $$Q = \chi(T),$$ where $\chi$ is equal to $1$ in a neighborhood of $t_1$  and is equal to zero in a neighborhood of $\sigma(T) \setminus \{t_1\}$. By Lemma \ref{idempotent}, it is sufficient to prove that  $Ran\; Q = L_1$.

With respect to the decomposition $H = L_1 \oplus M_1$ the operator $T$ is of the form $$T = \left(\begin{array}{cc} A & B \\0  & C \end{array} \right), $$
where $\sigma(A) = \{t_1\}\;$, $\sigma(C) = \sigma(T) \setminus  \{t_1\}$. Hence $$ Q = \chi(T) = \left(\begin{array}{cc} \chi(A) & \ast \\  0 & \chi(C) \end{array} \right) = \left(\begin{array}{cc} 1 & \ast \\ 0 & 0 \end{array} \right).$$ Hence $Ran\; Q = L_1$.
\end{proof}

\begin{lemma}\label{corner} Let $A \in B(H)$ be given by $A = (A_{ij})$ with respect to some orthonormal basis $\{e_i\}$ in $H$. If $|A_{11}| = \|A\|$, then $A_{1j} = A_{j1} = 0$, when $j\neq 1$.
\end{lemma}
\begin{proof} We have $$(AA^*)_{11} = \sum_{j\ge 1} |A_{1j}|^2 = \|A\|^2 + \sum_{j>1}|A_{1j}|^2.$$ Since $$ (AA^*)_{11} = (AA^*e_1, e_1) \le \|AA^*\| = \|A\|^2,$$ we get $A_{1j}$=0, when $j>1$. Applying this to $A^*$, we get $A_{j1}=0$, when $j>1$.
\end{proof}

\medskip

In what follows we will assume $|t_1|\ge |t_2|\ge \ldots$ which
always can be done by transposition of factors in the product.

\medskip

\begin{cor}\label{form} Suppose $T\in B(H)$ and   $$(T-t_N)^{k_N}(T-t_{N-1})^{k_{N-1}}\ldots (T-t_1)^{k_1}=0. $$ Then there exists $0\le m\le N$ such that with respect to the decomposition $H =  L_1 \oplus \ldots   \oplus L_m \oplus M_m$
$$ T = t_1 1 \oplus \ldots \oplus t_m 1 \oplus S,$$ where $S$ is such that $(S-t_{m+1})^{k_{m+1}} \ldots (S - t_N)^{k_N} = 0$ and $\|S\|>|t_i|$, for all $i\ge m+1$.
\end{cor}
\begin{proof} We write $T$ in upper-triangular form as in Lemma \ref{1} and then use Lemma \ref{corner}.
\end{proof}

\bigskip

Let $p$ be a polynomial in one variable, $C\ge 0$. Below the
universal $C^*$-algebra
 $$\mathcal A = C^* \left\langle p(x)=0, \; \|x\| \le C\right\rangle$$ is denoted by $\mathcal A$.

 \bigskip

\begin{theorem}\label{semiproj}  $\mathcal A$ is semiprojective.
\end{theorem}
\begin{proof} Write $p$ as $p(x) = (x-t_N)^{k_N}\ldots (x-t_1)^{k_1}=0$. Let $b\in A/I$, $I=\overline{\bigcup I_n}$, $(b-~t_N~)~^{k_N}\ldots (b-t_1)^{k_1}=0, \; \|b\| \le C$.
Embed $A/I$ into $B(H)$ and write $b$ as in Corollary \ref{form}.
Let $p_1, \ldots, p_m$ be the projections onto $L_1, \ldots, L_m$
and $p_{m+1}$ be the projection onto $M_m$. By Corollary
\ref{belong}, they all belong to $A/I$. Then  $$b = \sum_{i=1}^m t_i
p_i + s,$$ where $s\in  p_{m+1} A/I p_{m+1}$ satisfies the equation
\begin{equation}\label{polynom}(s-t_{m+1})^{k_{m+1}} \ldots (s - t_N)^{k_N} = 0\end{equation} and
\begin{equation}\label{radius}
\|s\|> \max_{i\ge m+1}|t_i|.\end{equation}
 By \cite{Blackadar} there
exists $n$ such that
   $p_i$'s can be lifted to projections $P_i$'s in $A/I_{n}$ with $\sum_{i=1}^{m+1}~P_i~=~1$.
   Since $A/I = (A/I_n)/(I/I_n)$, we have $$p_{m+1}\; A/I\; p_{m+1} = (P_{m+1}\; A/I_n\;
   P_{m+1})/(P_{m+1}\;
   I/I_n\;
   P_{m+1})$$
   and, by (\ref{polynom}), (\ref{radius}) and Lemma \ref{basic}
  we can lift $s$ to $S\in P_{m+1}\; A/I_n\;
   P_{m+1}$
 with $$(S-t_{m+1})^{k_{m+1}} \ldots (S - t_N)^{k_N} = 0,\;\; \|S\| = \|s\|  \le C.$$  Let \begin{equation}\label{triangular} a = \sum_{i=1}^m t_i P_i + S.\end{equation}
  It is a lift of $b$, $\|a\| \le C$ and $(a-t_N)^{k_N}\ldots (a-t_1)^{k_1}=0$. The last equality can be checked by direct
  calculations, but it is easier to say that  (\ref{triangular}) corresponds to upper-triangular form of $a$ as in Corollary \ref{form} and
  then the last equality follows instantly.
 \end{proof}

\begin{theorem}\label{RFD} $\mathcal A$ is RFD.
\end{theorem}
\begin{proof} Let $H = l^2(\mathbb N)$. We will identify the algebra $M_n$ of $n$-by-$n$ matrices  with $$B(l^2\{1, \ldots, n\})\subseteq
B(H).$$ Let $\mathcal B\subseteq \prod M_n$ be the $C^*$-algebra of
all $\ast$-strongly convergent sequences and let $\mathcal I$ be the
ideal of all sequences $\ast$-strongly convergent to zero. Then we
can identify $\mathcal B/\mathcal I$ with $B(H)$ by sending each
sequence to its $\ast$-strong limit.

We claim that any family $p_1, \ldots, p_n$ of projections  with
sum 1 in $B(H)$ lifts to a family of projections $P_1, \ldots, P_n$
with sum 1 in $\mathcal B$.  One way to prove this is to modify the
argument of Choi, used the proof of Theorem 7 of \cite{Choi}.  A more
modern approach is to use Hadwin's \cite{DonRFD} characterization of separable
RFD $C^*$-algebras:  in the unital separable case, $D$ is RFD if and only if
we can lift all elements of  $\hom_1(D, B(H))$ to $\hom_1(D, \mathcal B)$.
Clearly $\mathbb{C} ^n$ is RFD and this lifting problem for $\mathbb{C} ^n$
is equivalent to the needed lift of $n$ projections that sum to the identity.

Let $\pi: \mathcal A \to B(H)$ be the universal representation.
Arguments from the proof of Theorem \ref{semiproj} can be repeated
without any change to show that $\pi$ lifts to a $\ast$-homomorphism
$\tilde \pi: \mathcal A \to \mathcal B$. This lift gives a
separating family of finite-dimensional representations.
\end{proof}

\begin{theorem}\label{Calkin} Any $\ast$-homomorphism from $\mathcal A$ to Calkin
algebra lifts to a $\ast$-homomorphism to $B(H)$. In particular
$Ext(\mathcal A)=0$.
\end{theorem}
\begin{proof} As is well known, orthogonal projections (with sum 1) in Calkin algebra can be lifted
to orthogonal projections (with sum 1) in $B(H)$. Now we can repeat
arguments from the proof of Theorem \ref{semiproj}.
\end{proof}

\begin{remark}
To each $x$ in any $C^*$-algebra with $p(x)=0$ we can assign in a canonical
and functorial way a collection of projections that are orthogonal and sum
to one.  See Corollary~\ref{belong}, or \cite{Don}.  If $x \in A/I$ and these
projections lift, then $x$ lifts, preserving the relation $p(x)=0$ and
the norm.  In formal terms,
\[
\mathbb{C}^{N-1} \rightarrow C^* \left\langle x \;|\; p(x)=0, \; \|x\| \le C\right\rangle
\]
is conditionally projective.  Thus we have improved upon Theorem 2 in
\cite{Don} by incorporating the norm condition.
\end{remark}

\end{document}